\documentclass[12pt]{amsart}
\usepackage{epsfig}
\usepackage{verbatim}
\usepackage{amssymb}
\usepackage{amsfonts}
\usepackage{txfonts}
\newcommand{\ra}{{\rightarrow}}
\newcommand{\eproof}{\hfill\rule{2.2mm}{3.0mm}}

\newcommand{\Proof}{\noindent {\bf Proof.~~}}
\newcommand{\PG}{{\mathcal P}(G,q)}
\newcommand{\PGone}{{\mathcal P}(G_1,q_1)}
\newcommand{\PGtwo}{{\mathcal P}(G_2,q_2)}

\newcommand{\Z}{{\mathbb Z}}
\newcommand{\T}{{\mathbb T}}

\newcommand{\wgt}{{\rm wgt}}

\newcommand{\ao}{{\rm AO}}
\newcommand{\aof}{{\mathcal O}}
\newcommand{\aruso}{{\rm ARUSO}}
\newcommand{\dom}{{\rm dom}}

\renewcommand{\eqref}[1]{(\ref{#1})}
\newcommand{\eat}[1]{}

\newtheorem{prop}{Proposition}[section]
\newtheorem{lem}[prop]{Lemma}
\newtheorem{defi}{Definition}[section]
\newtheorem{coro}[prop]{Corollary}
\newtheorem{theo}[prop]{Theorem}
\newtheorem{Claim}[prop]{Claim}

\newtheorem{exam}{Example}[section]
\newtheorem{rmk}{Remark}[section]

\newcommand{\Nf}{\|f\|}
\newcommand{\A}{{\mathcal A}}

\def\A{\mathcal A}
\def\T{\mathcal T}
\def\>{>_{\sigma}}

\begin{document}
\baselineskip 18pt
\title{$G$-Parking Functions, Acyclic Orientations and Spanning Trees}
\author{Brian Benson}
\address{B. Benson, Department of Mathematics, University of Illinois, Urbana-Champaign, Illinois}
\email{benson9@illinois.edu}
\author{Deeparnab Chakrabarty}
\address{D. Chakrabarty, Department of Combinatorics \& Optimization, University of Waterloo,  Waterloo,
Ontario N2L 3G, Canada}
\email{deepc@math.uwaterloo.ca}
\author{Prasad Tetali}
\address{P. Tetali, School of Mathematics and School of Computer Science, Georgia Institute of Technology, Atlanta, GA 30332-0160.}
\email{tetali@math.gatech.edu}

\begin{abstract}
Given an undirected graph $G=(V,E)$, and a designated vertex $q\in V$,  the notion of a $G$-parking function (with respect to $q$) was independently developed and studied by various authors, and has  recently gained renewed attention.
This notion generalizes the classical notion of a parking function associated with
the complete graph. In this work, we study properties of 
{\em maximum} $G$-parking functions and provide  a new bijection between them and the set of spanning trees
of $G$  with no broken circuit. As a case study, we specialize some of our results to the graph corresponding to the discrete $n$-cube $Q_n$.
 We present the article in an expository self-contained form,  since we found  the combinatorial aspects of $G$-parking functions somewhat scattered in the literature, typically  treated in conjunction with sandpile models and closely related chip-firing games. 


\end{abstract}
\maketitle

\section{Introduction}
\setcounter{equation}{0}

The classical parking functions provide a bijective correspondence
between the spanning trees of the complete graph $K_n$ and certain integer-valued
functions on the vertices of $K_n$.  A notion of parking functions corresponding to
the spanning trees of an arbitrary graph $G$ is more recent and has been independently
developed in physics and combinatorics.  It was introduced by Bak, Tang and Wisenfeld
\cite{BTW87} as a self-organized sandpile model on grids, and was generalized to arbitrary graphs by D. Dhar \cite{Dhar90}.  See Definition~\ref{defn:park} below for the precise definition of a $G$-parking function, associated with a connected graph $G$. 

This notion is already rather powerful; besides generalizing the classical parking function from $K_n$ to an arbitrary graph, it has been investigated in the context of  chip-firing games \cite{Biggs99, Lopez97, PC} and the Tutte polynomial \cite{Biggs99b, CLe03} in discrete mathematics, and 
also investigated in algebra and related fields \cite{BN07, CRS02, PS04}.
However, some of the combinatorial aspects  of this topic appear somewhat scattered in the literature.

Several fundamental results concerning
the recurrent configurations of chip-firing can be derived without the chip-firing
context and terminology. For this reason, we shy away from introducing and discussing the chip-firing 
terminology.  Instead, in this article we describe various interpretations of  the $G$-parking functions  in the most elementary combinatorial ways.  Using a natural partial order $\prec$ on the set $\PG$ of parking  functions, we consider the maximal elements in this poset $\bigl(\PG, \prec\bigr)$ .  Much of our focus in this paper is on understanding the properties of such maximal parking functions.
The first result we describe (see Theorem~\ref{thm:main}) provides a new bijection between the maximal parking functions in the poset and  the set $\A(G;q)$ of acyclic orientations of $G$ with a unique source at $q$. 
En route, we describe what we call an Extended Dhar algorithm (since it is an extension of an algorithm due to Dhar \cite{Dhar90} to recognize $G$-parking functions) in providing an acyclic orientation corresponding to a maximal parking function.
We review  various combinatorial consequences  and algebraic connections of this correspondence.  For example, using known results (namely those of Greene and Zaslavsky \cite{GZ83} and more recent work of Gebhard and Sagan  \cite{GS99}),  we further identify a 1-1 correspondence between the set of maximal parking functions and the set of spanning trees with no  ``broken circuits," or equivalently, the set of ``safe" spanning trees ; see Section~\ref{sec:MPF} for the definitions of these terms.  In this paper, we provide a much simpler bijection (compared to \cite{GS99}) between the set of  safe trees and the set of acyclic orientations with a unique sink (or equivalently, a unique source). Furthermore we generalize this bijection to one between all spanning trees and all $G$-parking functions which {\em preserves} the bijection between safe trees and maximal $G$-parking functions.
We must remark here that other bijective proofs between the set of $G$-parking functions and the set of spanning trees of $G$ (for arbitrary connected $G$) have been given by Chebikin and Pylyavskyy \cite{CP05}. However, to our knowledge, the  simpler bijection we report here, in Theorem~\ref{main} below, and its generalization given in Theorem ~\ref{thm:bijgen},  are indeed new.

As an additional  contribution, we describe a simple way to generate maximal parking functions in the Cartesian product graph $G_1 \Box G_2$,  using maximal functions in the (factor) graphs $G_1$ and $G_2$. We then specialize our study to understanding the parking functions in the discrete $n$-cube $Q_n$ on $2^n$ vertices. By describing certain special constructions of maximal parking functions $f$ on $Q_n$, we obtain a natural description of a set, ${\rm dom}(f)$, of parking functions --  those dominated,  in the partial order given by $\prec$, by a special maximal parking function $f$.
Interestingly enough we shall deduce (see Theorem~\ref{thm:domnum}) that
\begin{equation}
\label{domnum}
|{\rm dom}(f)| = \prod_{k=2}^n k^{{n\choose k}}\,,
\end{equation}
while it is a well-known fact that
\begin{equation}
\label{totalnum}
 |{\mathcal P}(Q_n, q)| = \prod_{k=2}^n (2k)^{{n \choose k}} = 2^{2^n-n-1}\,\, \prod_{k=2}^n k^{{n \choose k}}\,.
 \end{equation}
 
Recall  that (\ref{totalnum}) corresponds  to the total number of spanning trees
of $Q_n$ (see equation (5.85) in \cite{StanVol2}), using the matrix-tree theorem and the explicit knowledge
of the corresponding eigenvalues, to help evaluate the determinantal formula.  In light of the fact that finding a bijective proof accounting for  the  number of spanning trees of $Q_n$ has been open for several years, we hope this is a nontrivial step towards such a proof.

The paper is organized as follows.
In Section~\ref{sec:Dhar}, we review some preliminaries, including  Dhar's burn criterion, which determines whether a given function is a parking function. In Section~\ref{sec:AO}, we show the bijection between maximum parking functions and acyclic orientations with a unique source.  In Section~\ref{sec:safe_trees}, we describe our new and simpler bijection between the set of acyclic orientations with a unique sink and the set of safe trees.
 In Section~\ref{sec:PROD}, we describe a construction of maximum parking functions on Cartesian products of graphs.  In Section~\ref{sec:CUBE}, we focus our study on the $n$-cube $Q_n$, and provide some explicit constructions of maximum parking functions and related bounds.  In Section~\ref{sec:diffuse}, for expository purposes we review a   bijection between {\em diffuse states} (introduced in the context of chip-firing) and acyclic orientations of a graph. 
We conclude with some remarks on research in future directions and a few open problems in Section~\ref{sec:CONCL}.

\section{$G$-Parking Functions and Dhar's Burn Criteria}
\label{sec:Dhar}
\setcounter{equation}{0}

In this section we recall the definition of a $G$-parking function and review
Dhar's (burning) algorithm that can be used to determine whether an integer-valued function on the vertices of $G$ is  a $G$-parking function.

\begin{defi}
\label{defn:park}
  For a connected graph $G$, a $G$-parking function relative to vertex $q \in G$ is a function $f:V(G) \ra \Z_{\geq -1}$ such that $f(q)=-1$ and for every non-empty $A \subseteq V(G) \backslash \{q\}$, there exists $v \in A$ such that $0 \leq f(v) < d_{\overline{A}} (v)$, where $d_{\overline{A}}(v)$ is the number of edges $e = vw$ with $w \notin A$.\end{defi}

\begin{rmk}  Herein, we have modified the standard definition of a $G$-parking function somewhat.  The function $f$ is now defined on all of $V(G)$ instead of restricted to simply $V(G) \backslash \{q\}$ in order to improve the compatibility between $G$-parking functions and Cartesian product graphs such as $Q_n$; due to this change, for $f$ to be a $G$-parking function, $f(q)=-1$ necessarily.
\end{rmk}

\begin{prop} If for a function $f:V(G) \backslash \{q\} \ra \Z_{\geq 0}$, for every non-empty {\em connected} subgraph $A \subseteq G \backslash \{q\}$,  there exists $v \in V(A)$ such that $f(v) < d_{\overline{A}} (v)$, then $f$ is a $G$-parking function.  \end{prop}

\Proof Assume that, for all connected $A \subseteq G \backslash \{q\}$, that there exists $v \in V(A)$ such that $f(v) < d_{\overline{A}}(v)$.  Proceeding by contradiction, suppose that there is some disconnected $B \subseteq G \backslash \{q\}$ such that $f(v) \geq d_{\overline{B}} (v)$ for every $v \in V(B)$. Consider then any connected component $C$ of $B$.  Since $C$ is connected we have, by the hypothesis of the proposition, that $f(v) < d_{\overline{C}}(v)$, for some vertex $v$ in $C$.  Thus $d_{\overline{B}}(v) < d_{\overline{C}}(v)$, implying that there is a vertex $u$ in $\overline{C} \backslash \overline{B}$ such that $v$ and $u$ are connected by an edge 
in $G$; otherwise, either $f(v) \geq d_{\overline{C}}(v)$ or $f(v) < d_{\overline{B}} (v)$. This contradicts the choice of $C$.\eproof
 
Throughout we assume that the reference vertex $q$ is fixed, and we always
consider parking functions with respect to this fixed vertex $q$, without
necessarily bringing  explicit reference to it.

A natural question to ask is whether a given integer-valued function on the vertices
of $G$  can easily be tested for being a $G$-parking function.  In the context of
the so-called sandpile models, Dhar \cite{Dhar90} provided an algorithm, which
can be interpreted as an efficient algorithm to test if a given function is a $G$-parking function. This was observed in \cite{CP05}, wherein the algorithm was reformulated
as follows. 
Let $f : V\setminus \{q\} \to \Z^+$.  We assume that $f(q)=-1$.

Step 1. Mark any unmarked vertex $v$ which has more marked neighbors than $f(v)$.

Step 2. Repeat Step 1 until no more vertices can be marked.

Step 3.  Declare $f$ to be a $G$-parking function if and only if all the vertices
have been marked.

We omit the proof of correctness of the algorithm (as an exercise),  which follows in a fairly straightforward way from the definition of a parking function.


On the set ${\mathcal P}(G,q)$ of parking functions on $G$ with respect to $q$,
there is a natural partial order we may define:
\begin{defi}
Given two parking functions $f,g \in \PG$, we say $g\prec f$, if $g(v) \le f(v)$, 
for all $v \in V$. The maximal elements in this partial order will be referred
to as {\em maximal}  parking functions. Finally, a parking function with
the largest sum is called a {\em maximum} parking function.

For  $f \in \PG$, let $\Nf := \sum_v f(v)$, and 
${\rm dom}(f) =  \{ g \in \PG : g \prec f\}$.  Furthermore, for  $f, g \in \PG$, 
let $f \wedge g$ be the function on vertices, defined as $f\wedge g(v) := \min \{f(v), g(v) \}$,  for each  $v \in V(G)$.
\end{defi}



While the following propositions are perhaps folk-lore, the first part of Proposition~\ref{prop:dom_park} and Proposition~\ref{prop:bb} appear as Lemmas~7 and 5 in \cite{Biggs99b}.

 \begin{prop}
 \label{prop:dom_park}
(a)  Let $f \in \PG$,  and suppose $g:V(G) \ra \Z_{\geq -1}$ such that $g(q)=-1$ and 
 $0 \leq g(v) \leq f(v)$, for $v \in G$ with $v \neq q$.  Then $g \in \PG$.

 (b)  If $f,g \in \PG$, then $f\wedge g \in \PG$ and   $\dom(f \wedge g)= \dom (f) \cap \dom(g)$.
\end{prop}

 \Proof  While Part (a) is clear from the definition of a $G$-parking function, it can also be
seen using Dhar's algorithm. 

 Note that $f \wedge g (q) = -1$ and so $f\wedge g \in \PG$, by Part (a). Now
 $f\wedge g \prec f$ and $f\wedge g \prec g$, so $h\in \dom(f \wedge g)$ implies
 that $h \in \dom(f) \cap \dom (g)$.  Finally consider, $h \in \dom(f) \cap \dom (g)$.  Then at each vertex $v \in V(G)$, $h(v) \leq f(v)$ and $h(v) \leq g(v)$, so $h(v) \leq \min \{f(v), g(v) \} = f \wedge g (v)$, proving that $h \in \dom (f \wedge g)$.

\eproof


\begin{prop}
\label{prop:bb}
For every connected graph $G=(V,E)$, every  $f\in \PG$, we have  $\Nf \le |E| - |V|$.
More over, the equality is always achieved.

\end{prop}

\Proof
We may prove this by induction on the number $n\ge 1$ of vertices $G$. The base case
consisting of $V=\{q\}$ and no edges is trivially true. For the induction case, let $n\ge 2$. Given $f\in \PG$, let $v_n$ be the last vertex to be marked by the Dhar algorithm.   Then consider $H:=G\setminus \{v_n\}$, the graph obtained by removing $v_n$ and its incident edges. $H$  is  connected, since Dhar successfully marks all vertices before $v_n$, and more over, the function $f_H$ denoting, $f$ restricted to $H$, is an $H$-parking function with respect to  $q\in H$. Thus 
we may apply the induction hypothesis to $H$ and $f_H$  and complete the
proof: 
\[\|f\| = f(v_n) + \|f_H\|  \le \sum_{v \neq v_n} f(v) \le 
d(v_n)-1 + \bigl(|E| - d(v_n)\bigr) - (|V|-1) \le |E| - |V| \,,\]
where we also used the fact that $f(v) \le d(v)-1$, for every $v$ and parking function $f$.

The proof also suggests that by assigning the maximum possible value, at each step
in Dhar's marking algorithm, one easily obtains a (maximum) parking function which achieves the upper bound.
\eproof

Note that the quantity ${\bf g}(G):=|E| - |V| + 1$ is sometimes referred to 
as the {\em cyclomatic number} or the {\em Betty number} of the graph 
and due to our convention of assigning $f(q):=-1$, we have ${\bf g}(G) - 1$, as the bound in the above proposition.\\

Proposition \ref{prop:dom_park} also gives us a simple (albeit 
not necessarily efficient), inclusion-exclusion method to relate the set of maximum parking functions with the set of all parking functions.  But this has to wait until
the next section, where we observe another basic fact concerning
the maximum parking functions.

\section{Maximal $G$-Parking Functions and Acyclic Orientations with a Unique Source}
\label{sec:AO}
\setcounter{equation}{0}

Given a graph $G$, the  notion of  an acyclic orientation of the edges of $G$ is classical, with an extensive literature.  The notion of an acyclic orientation with a
{\em unique source} at a fixed vertex is less well-studied.  
Let $\ao(G)$ denote the set of acyclic orientations of the graph $G$ and  let $\A(G,q)$ be the set of acyclic orientations of $G$ with a {\em unique source} at vertex $q$. Finally, let  ${\rm MP}(G,q)$ denote the set of maximum $G$-parking functions.


\begin{theo}
\label{thm:biject}
There exists a bijection between $\A(G,q)$ and the set ${\rm MP}(G,q)$ of maximum $G$-parking  functions.
\end{theo}
\Proof 
Given an acyclic orientation $\aof \in \A(G,q)$ with a unique source at $q$, define the function $f= f(\aof)$ on the vertices of $G$: let $f(v)$  be the indegree (in $\aof$)  of $v$ minus 1. We will show that this correspondence provides the necessary bijection.

In any orientation, the sum of the indegrees equals the number of edges.  
Hence $\sum_v f(v)  = |E| - |V|$.  By using Dhar's algorithm, we may show that $f$
is in fact a parking function: starting with $q$, we may repeatedly mark and remove the current set of source(s) in the acyclic orientation of the remaining graph; since a vertex
$v$ with value $f(v)$ becomes a source only  when all its $f(v)+1$ in-neighbors
have been marked and removed, the Dhar criterion is satisfied. Also observe that
the procedure stops only after marking all  the vertices, since every acyclic orientation
has at least one source.  Thus $f(\aof) \in {\rm MP}(G,q)$.

To see that $f(\aof_1) \neq f(\aof_2)$, whenever $\aof_1 \neq \aof_2$, 
simply recall that  an acyclic orientation is uniquely determined by its outdegree sequence: starting with the sinks, orient all edges into the sinks, remove the sinks, and repeat the process by subtracting one from the outdegrees of  the neighbors of the sinks.

The proof will be complete once we establish the onto property, that every maximum
parking function can be obtained this way.  Given a maximum parking function $f\in {\rm MP}(G,q)$, we will construct an  orientation $\aof(f)$ using the following modification of Dhar's algorithm, and will show that $\aof(f) \in \A(G;q)$, thus essentially providing an
inverse map to the above construction.

\noindent
{\bf The Extended Dhar Algorithm}.

{\em Input}: A maximal parking function $f \in {\rm MP}(G,q)$

{\em Output}: An acyclic orientation $\aof(f)$ with a unique source at $q$.

Step 1. Start with $v=q$. Orient all edges out of $q$.

Step 2. If there exists a vertex $v$ which accrued indegree($v$) equal to $f(v)+1$, mark $v$ and orient the remaining edges incident at $v$ outward from $v$.

Step 3. Repeat Step 2, until all vertices  are marked and all edges are oriented.

The correctness of the original Dhar algorithm guarantees that all vertices will eventually be marked -- indeed, the indegree($v$) equals the number of neighbors marked
before $v$; thus all edges will be oriented, meaning that $\aof(f)$ is an orientation of the edges of $G$.  Observe that the indegree
of a vertex $v$ equals $f(v)+1$.   Since $q$ is unique with $f(q):=-1$, it must be that $q$ is  the unique source.    It is also easy to see that $\aof(f)$ is acyclic --
if there were to be a cycle, considering the first vertex in the cycle which was
marked, we obtain a contradiction to the way the edges were oriented (in Step 2 above)
from  a  marked vertex. 
\eproof

\begin{rmk}
Upon completion of this work, we discovered (thanks to Matt Baker), that
Theorem~\ref{thm:biject} can also be derived using chip-firing games: As described in 
\cite{GR}, the notion of a so-called {\em diffuse state} (introduced by \cite{JS95} and see Definition~\ref{defn:diffuse} below)  helps relate chip-firing configurations to acyclic orientations. 
 Also thanks to an anonymous referee of an earlier version of this work \cite{BT08}, we learned that Theorem~\ref{thm:biject} is Lemma~10 (under the name of {\em allowable} orientations) of Biggs \cite{Biggs99b}, where it is mentioned that this in fact goes back to an even earlier result of Greene and Zaslavsky \cite{GZ83}.    
\end{rmk}

It is now easy to observe the following fact (which appears as Lemma~8 in \cite{Biggs99}).

\begin{coro}\label{cor:max_max}
Every maximal parking function is a maximum parking function.
\end{coro}

\proof
This follows from  the proof of correctness of the Extended Dhar algorithm described above -- if $f$ were maximal, but not maximum, then there must be a vertex
in Dhar's marking whose indegree is at least $f(v)+2$. But then we can increase
$f(v)$ by one, and obtain a valid parking function, contradicting the maximality of $f$.
\eproof

We now return to prove the simple result that was promised at the end of the previous
section.
\begin{coro}  Let $G$ be a finite graph with $k$ maximum $G$-parking functions.  Then there exist $G$-parking functions $f_1, \ldots, f_n$ such that the number of $G$-parking functions is $$\sum_{i=1}^n \pm | \dom (f_i) |\,,$$ where $n=2^k-1$ and the sign $+$ or $-$ is uniquely determined  by Proposition \ref{prop:dom_park}.
\end{coro}

\Proof  Let $g_1, g_2, \ldots , g_k \in \PG$ be all the maximum $G$-parking functions.
Trivially, every non-maximal parking function is dominated by some maximal parking function, and Corollary~\ref{cor:max_max} lets us observe that,
\[\PG = \cup_{i=1}^k \dom(g_i) \,.\]
Now we may simply use the inclusion-exclusion formula to count the size of the union of the above $k$ sets:
\begin{eqnarray*}
|\PG| & = & |\cup_{i=1}^k \dom(g_i)| \\
& = &\sum_i |\dom(g_i)| - \sum_{1\le i < j \le k} |\dom(g_i) \cap \dom(g_j)| + \cdots \\
& & + (-1)^k
|\dom(g_1) \cap \dom(g_2) \cap \cdots \cap \dom(g_k)|\\
& = & \sum_i |\dom(g_i)| - \sum_{1\le i < j \le k} |\dom(g_i \wedge g_j)| + \cdots + (-1)^k
|\dom(g_1 \wedge g_2 \wedge \cdots \wedge g_k)|\,,
\end{eqnarray*}
which, upon using Proposition~\ref{prop:dom_park}, completes the proof.\eproof


It would indeed be interesting to see if the above corollary can be used in 
making progress towards obtaining a bijective proof for the number of spanning trees of $Q_n$, the $n$-dimensional hypercube.  In Section \ref{sec:CUBE}, we take
a modest step towards it.

\begin{theo} 
\label{thm:unique}
Let $G$ be a simple, connected graph.  Then, for a fixed choice of $q$, $G$ has a unique maximum $G$-parking function if and only if $G$ is a tree.\end{theo}

{\bf Proof}.  If $G$ is a tree, then there is only one parking function with respect to
any $q$ since $G$ has no cycles and, thus, each vertex can have at most one marked neighbor in the Dhar algorithm.  (Note that this is in fact tautological if one uses the bijection between the parking functions and the spanning trees of $G$.) Hence there is only one maximum $G$-parking function.

The other direction is less obvious. However, observe that in light of the bijection established in Theorem~\ref{thm:biject}, it suffices to show the following. {\em Whenever $G$ is connected and contains a cycle,  then there are at least two acyclic orientations for $G$, with $q$ as the unique source.}  
This is easy to establish (for example, by considering the standard directed acyclic graph (DAG) representation of the graph), and we leave the proof as a simple exercise.
\eproof

See Remark~\ref{rmk:unique_fn}  below for another short (but indirect) proof of the above theorem.

\begin{coro} A simple, connected graph $G$ has a unique maximum parking function $f$ if and only if the range of $f$ is a subset of $\{-1,0\}$.\end{coro}

\Proof This follows from the observation that any tree has a unique parking function and, for any vertex $v \neq q$, $f(v)=0$.\eproof

\begin{coro} For every $G$, \ $ |\PG| = 1$ if and only if $|{\rm MP(G,q)}| = 1$.
\end{coro}

\section{Bijections of Maximal $G$-parking Functions}\label{sec:MPF}
We now augment  our Theorem~\ref{thm:biject} above to show a bijection between maximal $G$-parking functions
and spanning trees with no broken circuits.
Before stating the theorem, we need to recall several definitions from the literature.

The Tutte (or Tutte-Whitney) polynomial of a graph $G=(V,E)$ is the two-variable polynomial defined as 
 \begin{equation}
 \label{tutte_poly}
 {\mathcal T}_G(x,y) = \sum_{A\subseteq E} (x-1)^{\kappa(A)-\kappa(E)}\, 
 (y-1)^{|A| - n + \kappa(A)}\,,
 \end{equation}
 where $n=|V|$ and $\kappa(A)$ denotes the number of connected components of the graph on $V$ using  edgeset $A$.
 For $\lambda \in \Z^+$, the chromatic polynomial $\chi_G(\lambda)$ of a graph $G$ is defined  as the number of proper vertex colorings of $G$ using $\lambda$ colors.
For a general variable $\lambda$, the following relation between the chromatic polynomial and the Tutte polynomial
is well-known (see e.g. \cite{Bari77, Bari79,Tutte54, Welsh, W32}):
\begin{equation}
\label{eqn:chr_tutte}
\chi_G(\lambda)) =  (-1)^{n-\kappa(G)}\, \lambda^{\kappa(G)}\, 
{\mathcal T}_G(1-\lambda,0)\,,
\end{equation}
where $n=|V|$ is the number of vertices and $\kappa(G)$ is the number of connected components of $G$.

\begin{defi}
\label{broken_circuit}
Given a graph $G=(V,E)$ and an ordering of all the edges of $G$, 
 a broken circuit $B \subseteq E$ is any cycle (of edges) of $G$ minus the largest (according to the ordering) edge in the cycle.
 \end{defi}

Note that since every cycle contains (or gives rise to) a broken circuit,  a collection of edges {\em not} containing a broken circuit must necessarily be acyclic. 
Inspired by the terminology of Kenyon and Winkler \cite{KW}, we call a spanning tree $T$ {\em safe},  if it contains no broken circuits. That is, for all edges $e$ not in the tree, there is an edge in the unique cycle formed when $e$ is added to the tree, which is larger than $e$.
The notion of a broken circuit is  more general, and in fact explains the classical terminology: in the context of a matroid, an independent set of elements of the matroid obtained from a circuit, by removing the largest element (once again, according to some a priori global ordering of all the elements) of the circuit.  


\begin{theo}
\label{thm:main}
For every undirected, connected graph $G$, the following quantities are all the same.
\begin{itemize}

\item{(a)} The number of maximum $G$-parking functions with respect to $q$.

\item{(b)} The number of acyclic orientations with a unique source at $q$.

\item{(c)} The number of spanning trees with no broken circuits, or equivalently, with zero
external activity.

\item{(d)} The coefficient (up to sign) of the $\lambda$-term in the chromatic polynomial $\chi_G(\lambda)$.

\item{(e)} The value (up to sign)  ${\mathcal T}_G(1,0)$ of the Tutte polynomial ${\mathcal T}_G(x,y)$, evaluated at $x=1$ and $y=0$.

\end{itemize}
\end{theo}

\proof
Theorem \ref{thm:biject} shows the equivalence of (a) and (b).
The equivalence between (c) and (d) is classical and is part of Whitney's Broken Circuit
theorem \cite{W32}: that the chromatic polynomial on $n$ vertices is
given by 
\[ \chi_G(\lambda) = \sum_{r=0}^{n-1} (-1)^j m_r \lambda^{n-r},\]
where $m_j$ is the number of $r$-subsets of edges of $G$ which contain no broken circuit. The term $m_{n-1}$ corresponds to (the absolute value of) the coefficient
of $\lambda$; note that the  $n-1$-subsets under consideration being necessarily acyclic, correspond to spanning trees which do not contain a broken circuit. 

The equivalence between (b) and (d) is due to Greene and Zaslavsky \cite{GZ83}.
A (direct) bijective proof of the equivalence between (b) and (c) is given by Gebhard and Sagan \cite{GS99}, using a modification of an algorithm of Blass and Sagan \cite{BS86}. In Section \ref{sec:safe_trees} we provide a much shorter proof of the equivalence between (b) and (c).

The equivalence between (d) and (e) follows  from  (\ref{eqn:chr_tutte}), and using
$\kappa(G) =1$, for a connected $G$. 
The equivalence between (a) and (e) follows from results of \cite{Dhar90} and  \cite{Lopez97}, which confirmed a conjecture
of Biggs \cite{Biggs99} in the context of chip-firing.  An inductive proof (using edge deletions and contractions) without involving chip-firing is due to Plautz and Calderer \cite{PC}. As described in \cite{PC}, the work of Dhar and Lopez provides the following
result:
\[ {\mathcal T}_G(1,y) = \sum_{f\in \PG} y^{w(f)}\,,\]
where $w(f) = |E(G)| - |V(G)| + \|f\|$, hence the equivalence of (a) and (e).
The results in \cite{Lopez97} and \cite{CLe03} also establish the equivalence between (c) and (e), with the minor modification that broken circuits are equivalently described using
external acitivities - each broken circuit contributes an external activity of one to a spanning tree.
\eproof

In addition to the above,  Cori and Le Borgne \cite{CLe03}  describe certain decreasing traversals of vertices and edges and  a notion of strong edges to provide a bijection
between recurrent chip-firing configurations (with a fixed ``level")
and spanning trees with a fixed ``external activity".  While the level corresponds to the
sum of the values of a parking function (up to an additive shift), the external activity
reflects the number of broken circuits, and we refer the interested reader to
their paper for additional information.

\begin{rmk}
\label{rmk:unique_fn}
Observe that the nontrivial part of Theorem~\ref{thm:unique}
follows easily using the equivalence between (a) and (c).
Indeed, let $G$ be  a connected graph which contains a cycle of length 3 or more. 
Then given a spanning tree which contains no broken circuit, we include
an edge not in the tree to form a cycle $C$. Since $C$ is of length 3 or more, there
must be an  edge, which is not the largest edge in the cycle, that can be removed, 
giving another spanning tree with no broken circuits. Hence a (connected) graph containing a cycle has more than one parking function.
\end{rmk}


\subsection{Bijection between Acyclic Orientations with Unique Sink and Safe Trees} \label{sec:safe_trees}

In this section we give a shorter proof, of  equivalence of (b) and (c) of Theorem ~\ref{thm:main}, than the ones reported in Gebhard and Sagan \cite{GS99} and Gioan and Las Vergnas \cite{GV05}. 
Note that the bijection in \cite{GS99} is not
activity-preserving while the one in \cite{GV05} is; also see \cite{LV84}.

For ease of presentation, we will consider orientations with a unique sink rather than a unique source -- clearly this is equivalent. 

Let $\sigma$ be any total ordering of the edges of $G$. Given two edges $e$ and $f$, we say $e$ is {\em larger} than $f$ if $\sigma(e) > \sigma(f)$.
Similarly we say $e$ is {\em smaller} than $f$, if the inequality is otherwise. Recall, a spanning tree $T$ of $G$ is called {\em safe} with respect to $\sigma$ if for any edge  $e\notin T$, there exists at least one edge $f$ in the unique cycle in $T+e$ such that $f$ is larger than $e$. Let $\T(G,\sigma)$ be the set of safe trees with respect to $\sigma$, and let $\A(G,q)$ be the set of acyclic orientations of $G$ with $q$ being the {\em unique sink}.

\begin{theo}\label{main}
For any total order $\sigma$, there exists a bijection $\mu_\sigma: \T(G,\sigma) \to \A(G,q)$\,.
\end{theo}
\noindent
Henceforth, we fix $\sigma$ and do not write it as a subscript.
Before we proceed, we make a few more definitions and observations. An arborescence with root $q$ is a directed spanning tree with all vertices except the root having out-degree exactly $1$, and the root having out-degree $0$. Any spanning tree corresponds to a unique arborescence with root $q$ and henceforth we will use the terms interchangeably. Given a vertex $i$, we let $P_i$ denote the unique directed path from $i$ to $q$. Given vertex $i$ and $j$, we let $meet(i,j)$ be the first vertex in the intersection of $P_i$ and $P_j$. That is, the path from $i$ to $meet(i,j)$ and the path from $j$ to $meet(i,j)$ are disjoint except at $meet(i,j)$. If $j$ lies on $P_i$, we let $j$ be $meet(i,j)$. 
Observe that for three vertices $i,j,k$, either (a) $meet(i,k) = meet(j,k)$, or (b) $meet(i,k)$ lies on $P_i$ and $meet(j,k) = meet(i,j)$, or (c) $meet(j,k)$ lies on $P_j$ and $meet(i,k) = meet(i,j)$. Given two vertices $i$ and $j$, we will denote the largest edge in the path from $i$ to $meet(i,j)$ as
$e_{ij}$ and the largest edge from $j$ to $meet(i,j)$ as $e_{ji}$. If $j = meet(i,j)$, we let $e_{ji}$ be the null edge.

We (abuse notation and) say $i\>j$  if $e_{ij}$ is larger than $e_{ji}$. We will also define $e_{ii}$ to be a null edge.
Note that if $i\>j$ and $j\>k$, then going over the three possibilities of $meet(i,k)$ we see that $i\>k$. Thus $\>$ is transitive and induces a total ordering of vertices. We say $i$ {\em dominates} $j$, if $i\>j$. It is instructive to note that $i$ dominates all vertices in $P_i$.

We now describe a mapping $\mu$ from all arborescences with root $q$ to acyclic orientations with unique sink $q$.
We will prove that distinct safe trees lead to distinct arborescences. This proves that the mapping $\mu$ restricted to safe trees is one-to-one.
Furthermore, given an acyclic orientation, we describe a procedure $\pi$ which takes an acyclic orientation and returns a safe tree. Moreover,
for any orientation $\mathcal{O} \in \A(G,q)$, we have $\mu(\pi(\mathcal{O})) = \mathcal{O}$. This shows  that $\mu$ is onto and thus it is a bijection.\\
\\
\noindent
{\bf Arborescence to Acyclic Orientation ($\mu$)}: Given an arborescence, orient an edge $(i,j)$ as 
$i$ to $j$ if $i\> j$, or vice-versa. By the transitivity of $\>$ it is clear that the orientation is acyclic. Also every vertex dominates the root which therefore is the unique sink. We will call an orientation so obtained as one {\em induced} by the arborescence.  The following lemma will show that two safe arborescences cannot lead to the same orientation.
\def\O{\mathcal O}
\begin{lem} 
Let $T_1$ and $T_2$ be two distinct safe arborescences and $\O_1 = \mu(T_1)$ and $\O_2 = \mu(T_2)$. Then $\O_1 \neq \O_2$.
\end{lem}
\begin{proof}
We prove the contrapositive: suppose $\O_1 = \O_2 = \O$, then
we show that $T_1 = T_2$. Consider the trees rooted at $q$ (note that $q$ is the unique sink of $\O$)
with edges directed towards $q$. We now show that  for each vertex of the graph:

\noindent
(*) The unique out-neighbor in
$T_1$ is the same as that in $T_2$ which will imply that both trees are the same.

\noindent
Since $\O$ is acyclic with a unique sink, the vertices $V$ can be decomposed as $V= ({q} =: S_0 \cup S_1 \cup ... \cup S_r)$
for some $r\ge 1$, where $S_i$ is the set of vertices which are sinks in the digraph $G\setminus \bigcup_{\ell=0}^{i-1} S_\ell$.
Let $S_i$ be the first set (with the least $i$) to contain a vertex violating (*).
Let this vertex be denoted $i$ (abusing notation). Let $(i,j)$ and $(i,k)$
be the unique out-neighbors of $i$ in $T_1$ and $T_2$ respectively,
with $j\neq k$, $j\in S_j$ and $k\in S_k$ (again abusing notation).
Observe that $j, k < i$ (that is, $S_j,S_k$ precede $S_i$), by the nature of the decomposition,  and the fact that $\O=\mu(T_1)$ and $\O=\mu(T_2)$.

Now consider the undirected cycle using the edges $(i,j), (i,k)$ and the unique
paths, $P_j$ from $j$ to $q$ and $P_k$ from $k$ to $q$. Observe that
both $P_j$ and $P_k$ are contained in $T_1$ {\em and} $T_2$, by the choice of $i$.
Also note that the largest edge in this cycle
must be either on $P_j$ or on $P_k$, since both $T_1$ and $T_2$ are safe!
We now get a contradiction -- if the largest edge is on $P_j$, then our definition
of $\mu$ demands that $(i,j)$ be oriented from $j$ to $i$ in $T_1$; similarly if it is on $P_k$, then the edge $(i,k)$ be oriented from $k$ to $i$ in $T_2$.
\end{proof}

\noindent
{\bf Acyclic Orientations to Safe Trees($\pi$): } 
Now we describe a procedure to get a safe tree from an acyclic orientation $\O$ with single sink $q$. 
Let $d(v)$ denote the out-degree of vertex $v$. Note that $d(q) = 0$. For an edge oriented $i$ to $j$, we say $j$ is a out-neighbor of $i$
and $i$ is an in-neighbor of $j$. At each step we maintain a set of {\em labeled} vertices $X$ and an arborescence $T$ spanning $X$.
We maintain the invariant that there is no edge from a vertex in $X$ to a vertex in $Y := V\setminus X$, (think of  $X$ as a large sink). Initially, $X = \{q\}$.
Note that since the orientation is acyclic, at any step there is at least one vertex $u$ in $Y$ which is a sink in the induced graph $G[Y]$, that is, $u$ has no out-neighbors in $Y$. We pick one such $u$ arbitrarily.  Also, since there is a unique sink $q$, this vertex $u$ must have at least one out-neighbor in $X$. We add $u$ to $X$ and we connect $u$ to the arborescence $T$ as follows. 

Let $X_u \subseteq X$ be the set of out-neighbors of $u$. Let $x\in X_u$ be the vertex which dominates all other vertices in $X_u$ with respect to the current arborescence $T$. Let $W_u \subseteq X_u$ be the subset of all vertices $v$ such that $(u,v) \> e_{xv}$. That is, the edge $(u,v)$ is larger than the largest edge in the path from $x$ to $meet(x,v)$. Note that $x\in W_u$ and thus $W_u$ is non-empty. Connect $u$ to the vertex $v$ with the largest $(u,v)$ among all $v$ in $W_u$.  We end when $X = V$ with an arborescence $T$.

\begin{lem} 
The arborescence $T$ obtained at the end is a safe tree. Moreover, $\mu(\pi(\O)) = \O$, for $O\in \A(G, q)$.
\end{lem}
\begin{proof}
Consider an edge $(j,i)$ not in the tree. Note that at each step exactly one node is added to the arborescence. Also note that for an arc oriented $(j,i)$ in $\O$,  $i$ is added before $j$. Let $X_j$ be the set of labeled vertices in the step when $j$ is added to the arborescence. Note that $i\in X_j$. Suppose $x\in X_j$ was the dominator of $X_j$ and $(j,k)$ was the edge added at this step.

Observe that whenever the procedure adds a new vertex $j$,
the vertex $j$ dominates all  other vertices in the arborescence. This is because $(j,k) \> e_{xk}$, that is, the largest edge in the path $x$ to $meet(x,k) = meet(x,j)$ and thus $j$ dominates $x$ and so every other vertex. Thus the tree $T$ induces the same orientation $\O$. It remains to show that $T$ is safe, that is, $(j,i)$ is not the largest edge in the cycle $T+(j,i)$.

If $i\in W_j$, then by choice of $k$ (based on how the edge $(j,k)$ was added above), $(j,k) \> (j,i)$ and thus $(j,i)$ is not the largest cycle in $T+(j,i)$. 

If $i\notin W_j$, this means $(j,i) < e_{xi}$. That is, $(j,i)$ is smaller than the largest edge in the path from $x$ to $meet(x,i)$. Also, by definition, 
$(j,k) > e_{xj}$ that is, $(j,k)$ is larger than the largest edge in the path from $x$ to $meet(x,k)$.  Let the path from $x$ to $meet(x,i)$ be $Q_i$ and that from $x$ to $meet(x,k)$ be $Q_k$.

Now consider $meet(i,k)$. If $meet(i,k)$ does not lie on $P_x$, then $meet(x,i) = meet(x,k)$. Thus, $(j,k) \> (j,i)$. 
If $meet(i,k)$ lies on $P_x$, then it must be $meet(x,i)$ or $meet(x,k)$. If the latter, then $Q_i \subseteq Q_k$, and therefore $(j,k) \> (j,i)$.
If $meet(i,k) = meet(x,i)$, then $Q_k \subseteq Q_i$. Either the largest edge in $Q_i$ lies in  $Q_k$ and we are done as before;
Or, the largest lies in $Q_i\setminus Q_k$, which lies in the cycle formed in $T+(j,i)$. Thus, $(j,i)$ is smaller than the largest edge in the cycle of $T+(j,i)$, completing the proof. \end{proof}

\begin{rmk}
\label{new_remark}
Thanks to an anonymous referee, we learnt that Biggs and Winkler \cite{BW97} had actually  given a simple bijection somewhat similar to the one we describe above, but unlike ours, their bijection does not preserve the so-called external activity. This is an important distinction for us, since we extend the above notion in the next subsection to provide a more general bijection.
\end{rmk}

\subsection{Extension to a bijection of $G$-Parking Functions and Spanning Trees}\label{sec:EXT}
\def\P{\mathcal P}
In this section we generalize the above 
to a bijection between $G$-parking functions and all spanning trees of $G$, in such a way that it  preserves the bijection between maximal parking 
functions and safe trees. We use the definitions of the previous section and make a few more definitions and claim below before demonstrating the bijection.

Given an arborescence $T$ spanning only a subset of vertices $X\subseteq V$, and a vertex $u\notin X$, we make a few definitions and observations
which will be useful in our bijection. Let $\Gamma(u)$ be the neighbors of $u$ and let $X_u:= \Gamma(u) \cap X$. We now describe an order on 
the vertices $X_u$ (which could be different from the total order of the previous paragraph, but is related), which we call the {\em power order} of $X_u$.
Intuitively, given two neighbors $v$ and $w$ of $u$ in $X_u$, if $v$ is more powerful than $w$, then in the tree $T+(u,v)$, $u$ would dominate $w$
{\em and} the edge $(u,w)$ doesn't form a broken circuit with $T+(u,v)$.  Moreover, in the tree $T+(u,w)$, either $v$ dominates $u$ {\em or} the edge $(u,v)$ forms a broken circuit with $T+(u,w)$.

\begin{figure}[h]
\centering
\includegraphics[scale=0.5]{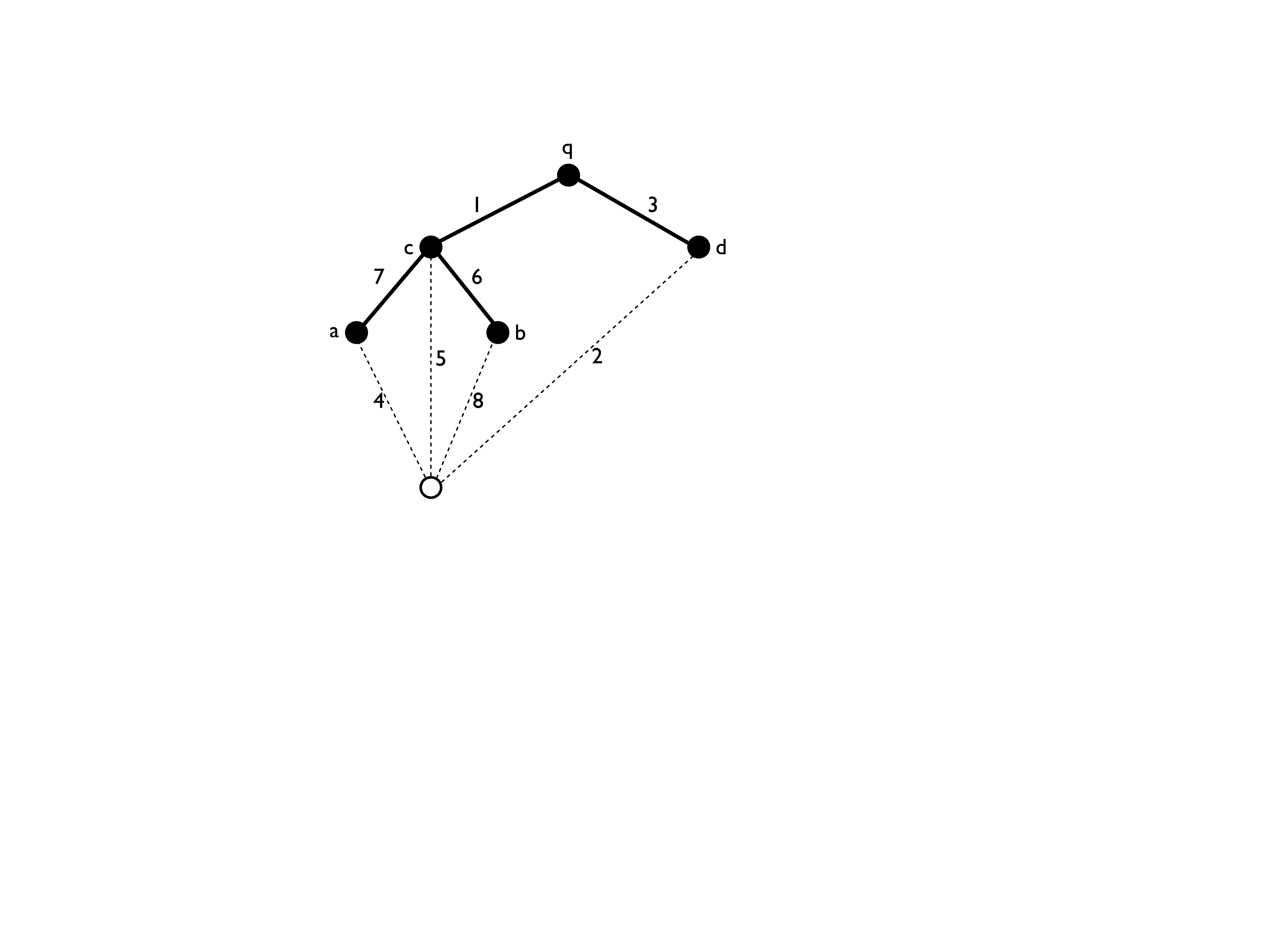}
\caption{The numbers on the edges correspond to the $\sigma$ value. The dark edges form the tree $T$ and $u\notin T$. Note that $a\> b\>d\>c$ in $T$. However, the power order of $X_u = \{a,b,c,d\}$ is as follows. Firstly, $x=a$ and $W_u = \{a,b\}$. Since $(u,b)$ is bigger, $b$ is the most powerful. Second powerful is $a$. After that $x=d$ and $W_u = \{c,d\}$, and thus the complete power order is $(b,a,c,d)$. }
\label{fig:powerorder}
\end{figure}

Let $x$ be the vertex in $X_u$ which dominates all other vertices in $X_u$ with respect to the current arborescence $T$. 
Let $W_u \subseteq X_u$ be the subset of all vertices $v$ such that $(u,v) \> e_{xv}$. That is, the edge $(u,v)$ is larger than the largest edge in the 	path from $x$ to $meet(x,v)$.  Note that $x\in W_u$ and thus $W_u$ is non-empty.  Let $v$ be such that $(u,v)$ is largest among all $v$ in $W_u$. 
Call $v$ the most {\em powerful} element of $X_u$. Delete $v$ from $X_u$ and repeat till one gets an order on all vertices of $X_u$. We call this order
the {\em power order w.r.t $T$}. In Figure \ref{fig:powerorder} we give an illustrative example.

\begin{Claim}\label{claim:power}
Given a tree $T$ and a vertex $u\notin T$ with neighbors $X_u$ in $T$, and two vertices $v$ and $w$ in $X_u$.
If $v$ is more powerful than $w$, then $u$ dominates $w$ in the tree $T+(u,v)$ and $(u,w)$ doesn't form a broken circuit
with $T+(u,v)$. 
 In the tree $T+(u,w)$, either $v$ dominates $u$ or $(u,v)$ forms a broken circuit with $T+(u,w)$.
\end{Claim}
\begin{proof}
Consider the time when the power order of $v$ is determined. Let $x$ be the dominator at that stage and let $X_u$ the set of neighbors
of $u$ remaining. Note $w\in X_u$. Either $v\in W_u$ or $v=x$ and $W_u$ is empty. 
If the former, $(u,v)$ is larger than $e_{xv}$ and thus in $T+(u,v)$, $u$ dominates $x$ which dominates $w$. Moreover $(u,w)$ is either smaller than $(u,v)$ if $w\in W_u$, or smaller than $e_{xw}$ which is smaller than $e_{uw}$ since $u$ dominates $x$. In any case, $(u,w)$ doesn't form a broken circuit in $T+(u,v)$. Also, in $T+(u,w)$, $(u,v)$ forms a broken circuit since it is larger than the largest edge from $v$ to $w$.

If the latter, then $u$ dominates $w$ in $T+(u,v)$ since $v=x$ dominates $w$ in $T$. Also since $W_u$ is empty, the edge $(u,w)$ is smaller than $e_{vw}$
implying $(u,w)$ doesn't form a broken circuit in $T+(u,v)$. Also, in $T+(u,w)$, $v$ dominates $u$ since $(u,w)$ is smaller than $e_{vw}$.
\end{proof}

Now we are ready to present the next theorem which gives the desired bijection between $G$-parking functions and spanning trees of $G$.
Let $\T(G)$ be the set of all spanning trees of $G$. Recall that $\T(G,\sigma)$ was the set of safe spanning trees with respect to $\sigma$.

\begin{theo}\label{thm:bijgen}
There exists a bijection $\pi: \P(G,q) \to \T(G)$
such that for all $f\in \P(G,q)$ which is maximal, $\pi(f)$ is a safe tree.
\end{theo}

As in the proof of Theorem \ref{thm:main}, we describe mappings $\mu:\T(G) \to \P(G,q)$ and $\pi: \P(G,q) \to \T(G)$ and show that for any parking function $f$, $\mu(\pi(f)) = f$ (Lemma \ref{lem1}) and for any spanning tree $T$, $\pi(\mu(T)) = T$ (Lemma \ref{lem2} ). Furthermore we show that for a tree $T$, $\mu(T)$ is a maximal parking function iff $T$ is safe (Proposition \ref{prop:safemax}). This completes the proof of Theorem \ref{thm:bijgen}\\

\noindent
{\bf Spanning Trees to Parking Functions ($\mu$):}
 \noindent
 Given a spanning tree $T$, let $E'$ be the subset of edges not in $T$ 
which form a broken circuit with $T$. Delete $E'$ from $E$.
Direct the edges of $T$ with $q$ as the unique source -- that is,
all vertices except $q$ have an out-degree of $1$ and $q$ has an out-degree of $0$. Given this arborescence, 
for any undirected edge $(i,j)$ in $E\setminus E'$, orient it as $(i,j)$ if $i\> j$, or $(j,i)$ if $j\> i$. This gives di-graph $D$.
The parking function  $f := \mu(T)$ is defined as $f(v) = \mbox{out-degree}_D(v) - 1$ in the orientation of the edges of $E\setminus E'$.
\begin{prop}\label{prop:safemax}
$f:=\mu(T)$ is a maximum parking function iff $T$ is a safe tree.
\end{prop}
\begin{proof}
Note that if $T$ is safe, $E'$ is empty and therefore all the edges are oriented. Thus $\sum_v f(v) = |E| - |V|$ and by
Proposition \ref{prop:bb}, $f$ is a maximal parking function. On the other hand if $T$ is not safe, then $E'$ is not 
empty and thus  $\sum_v f(v) < |E| - |V|$ implying $f$ is not a maximal parking function.
\end{proof}
\noindent
{\bf Parking Functions to Spanning Trees ($\pi$):}
\noindent
Given a valid parking function, we use a modification of Dhar's algorithm to obtain the spanning tree.
We maintain a set of vertices $X_i$ connected via an arborescence $T_i$, with $X_0$ initialized to $\{q\}$ and $T_0$ is $\varnothing$. 
At each step we add one vertex to $X_i$ and one edge to $T_i$. In the end we get $X_{n-1} = V$ and $T_{n-1}$ is the spanning tree returned. 
We describe the $i+1$th step.
Let $\Gamma(v)$ denote the neighbors of $v$.
\begin{enumerate}
\item Let $S_{i+1} := \{v\in V\setminus X_i: |\Gamma(v) \cap X_i| > f(v)\}$. 
\item 
         For every $u\in S_{i+1}$, $X_u := \Gamma(u) \cap X_i$ and let $M(u)$ be the $(|X_u| - f(u))$th vertex in the power order of $X_u$ with respect to
         $T$. Let $Y_u$ be the vertices in $X_u$ more powerful than $M(u)$ and $E_u$ be the set of $|X_u| - f(u) - 1$ edges of the form $(u,v)$ where 
         $v\in Y_u$.      
\item Note that adding all edges of the form $(u,M(u))$ to $T_i$ gives a new tree $T'$. In $T'$, let $u$ be the vertex in $S_{i+1}$
which is dominated by all other vertices in $S_{i+1}$ with respect to $T'$.  Add $u$ to get $X_{i+1}$ and the edge $(u,M(u))$ to get $T_{i+1}$.
\end{enumerate}

\begin{Claim}\label{claim:order}
Given a parking function $f$, let the vertices be added in order 
$\{q = u_0,u_1,\cdots,u_{n-1}\}$. That is $X_i := \{u_0,\cdots,u_i\}$.
Then with respect to the  tree $T_i$, $u_i$ dominates $u_{i-1}$ which dominates
$u_{i-2}$ and so on. 
\end{Claim}
\begin{proof}
The proof is by induction. At stage $i$,
let $S_i$ be the set of vertices as defined above. Suppose $u_i$ is 
added at this stage. We will be done if we show $u_i$ dominates $u_{i-1}$.

Two cases arise: If $u_{i-1}$ is not a neighbor of $u_i$, then $u_i$ must have
been in $S_{i-1}$ as well implying $u_i$ dominated $u_{i-1}$ in $T'_{i-1}$ and since
$M(u_i)$ doesn't change, dominates $u_{i-1}$ in $T_i$ as well.

If $u_{i-1}$ is a neighbor of $u_i$, then either $|X_{u_i}| = f(u_i) + 1$ and $u_i$ 
connects to the most powerful of its neighbors in $X_{u_i}$ and thus dominates
$u_{i-1}$ which is in $X_{u_i}$. Or, $|X_{u_i}| > f(u_i) + 1$, which once again implies
$u_i$ was in $S_{i-1}$ and moreover, either $M(u_i)$ becomes $u_{i-1}$, or still remains 
more powerful than $u_{i-1}$. In any case, $u_i$ dominates $u_{i-1}$ in $T_i$.
\end{proof}	

\begin{Claim}\label{claim:mainpi}
For any vertex $u$, the set of edges $E_u$ are precisely the set of edges from $u$ to $X_u$ which form broken circuits with $T$. 
\end{Claim}
\begin{proof}
Consider vertex $u$ connecting to vertex $w$ in $X_u$. From the first part of Claim \ref{claim:power} we see that the edges not in 
$E_u$ do not form broken circuits. Let $(u,v)$ be an edge in $E_u$. Since $v$ is more powerful than $w$, from the second part of Claim \ref{claim:power}
we see either $(u,v)$ forms a broken circuit with $T$ or $v$ dominates $u$ in $T$. The second possibility is precluded by Claim \ref{claim:order}.
\end{proof}

\begin{lem}\label{lem1}
For any parking function $f$, we have $\mu(\pi(f)) = f$.
\end{lem}
\begin{proof}
Let $T$ be the arborescence formed by rooting the tree $\pi(f)$ at $q$. Fix a vertex $u$. From Claim \ref{claim:mainpi} we have that
$\mu$ will first remove all the edges in $E_u$. After the removal of these edges, 
$u$ will dominate the remaining $f(u) + 1$ vertices in $X_u$ and thus its out-degree will be that. Thus, $\mu(\pi(f))(u) = \mbox{out-degree}_D(u) - 1 = f(u)$.
\end{proof}

\begin{lem}\label{lem2}
For any spanning tree $Z$, we have $\pi(\mu(Z)) = Z$.
\end{lem}
\begin{proof}
Let $f := \mu(Z)$. Abuse notation and call the arborescence obtained by rooting $Z$ at $q$, also $Z$. 
Let $E'$ be the set of edges which form broken circuits with $Z$. Let $D$ be the di-graph obtained 
by orienting the edges of $E\setminus E'$ with respect to $Z$. Given a subset of vertices $X$, let $Z[X]$ be the induced 
sub-forest of $X$.
 
The proof proceeds by induction on the stages of the algorithm computing $\pi$. We assume at stage $i$, the current tree of the algorithm, $T$,
is a subtree of the tree $Z$. That is $T = Z[X_i]$. We also assume that for every vertex in $X_i$, all its out-neighbors in $D$ are also in $X_i$.
These are vacuously true at stage $0$. We now show that at stage $i+1$ a vertex connects to tree $T$ using an edge of $Z$ and all its out-neighbors
of $D$ are in $X_i$.

Consider vertices in $V\setminus X_i$. At least one of these vertices $v$ must have all its $f(v)+1$ out-neighbors of $D$ in $X_i$ for otherwise
we would get a cycle in $D$. Call this set of vertices $S$. Note that $Z[X_i \cup S]$ is connected as the vertices in $X$ can only connect to $Z$ using some edge of $D$.
Let $u^*\in S$ be the vertex which is dominated by all other vertices in $S$ in $Z[X_i \cup S]$. Let $(u^*,w^*)$ be the edge in $Z$ with $w^*\in X_i$.
We claim that the algorithm which computes $\pi$ also picks $u^*$ in this stage and $w^*$ is $M(u)$.

Note that $S\subseteq S_{i+1}$ since each vertex in $S$ has at least $f(u) + 1$ neighbors in $X_i$. In fact, we show for every vertex $u\in S$, 
the corresponding $w$ in $X_i$ where $(u,w)\in Z$ is in fact $M(u)$. Call the set of these $f(u)+1$ out-neighbors of $u$ in $X_i$, $Y_u$.
Observe that any edge of the form $(u,v)$ with $v\notin Y_u$ must form a broken cycle with $Z$. This is because these edges
are not directed towards $v$ and cannot be directed towards $u$ by the induction hypothesis. Since they form broken 
cycles with $Z$, these $|X_u| - |Y_u|$ vertices must be more powerful than $w$ with respect to the arborescence $Z$ restricted to vertices
of $X_i$, that is $T$ by the induction hypothesis. This follows from the definition of power. Moreover, $w$ must be powerful than all other vertices  
of $Y_u$ with respect to $T$ since $u$ dominates all these vertices. Thus $w$ is the $(|X_u| - |Y_u| +1)$th powerful vertex in $X_u$ with respect to $T$, that is,
$w = M(u)$.

We will be done if we show any vertex $u'\in S_{i+1}\setminus Q$ dominates some vertex in $S$ with respect to $Z$. If this is the case,
then the algorithm would choose the vertex which is dominated by all vertices in $S_{i+1}$ and it has to be the vertex $u^*$.
But this is true since $u'$ has some out-neighbor of $D$ in $V\setminus X_i$ --  a path following argument shows we must reach a vertex
$v\in S$ from $u'$ using edges in $D$. In other words, $u'$ dominates $v$ with respect
to $Z$.
\end{proof}

In Figure \ref{fig:example} below, we give an example of the bijection on a simple $4$-vertex $5$-edge graph.
\begin{figure}[h]
\includegraphics[scale=0.5, angle=90]{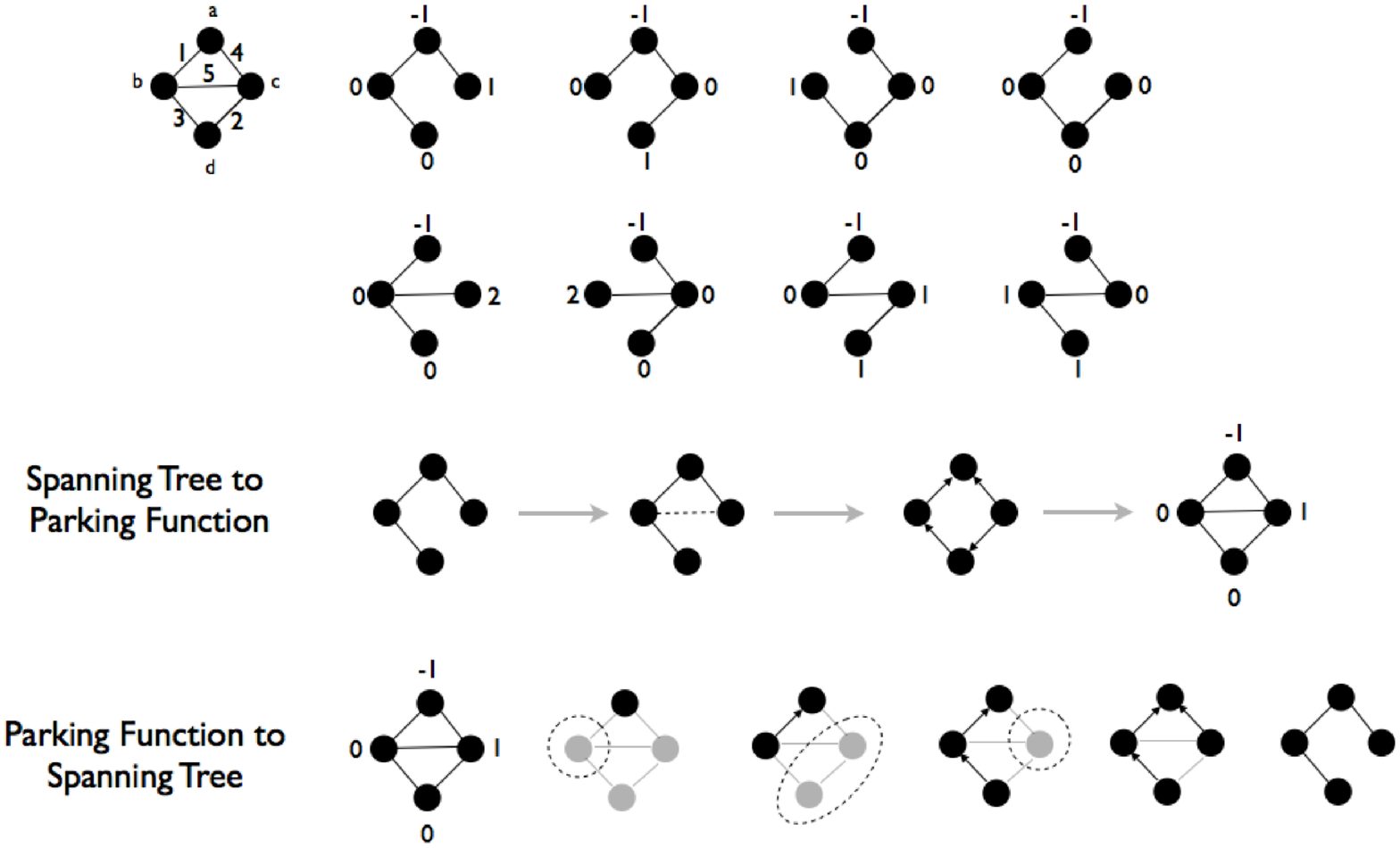}
\caption{We illustrate our bijection on the graph on the top left corner. The numbers on the edges is the ordering $\sigma$. The graph has 
$8$ spanning trees and they are shown beside the graph. The numbers on each spanning tree correspond to the $8$ possible parking functions
of the graph. It is instructive to note that the second row of spanning trees correspond to safe trees w.r.t to $\sigma$ and correspond to maximum parking 
functions -- note the sum of the numbers is $1$ for the bottom $4$ trees. We now take the first tree and show how it leads to the parking function 
via the function $\mu$ and then how vice-versa is obtained via the function $\pi$. 
$\mu$ first recognizes the edges which form broken circuits with $T$ -- the dotted edge is the only one in this case. Then it orients the tree edges towards the root $a$. It orients the edge $(c,d)$ towards $d$ because $c\>d$ w.r.t $T$. After the orientation, the parking function is found by subtracting $1$ from the out-degrees. 
Going from the parking function to the tree, the dark vertices denote the set $X_i$ at each step, while the dotted lines enclose the set $S_i$.
When $S_i$ has size more than $1$ (when it contains vertices $\{c,d\}$), it finds $M(c) = b$ and $M(d) = b$ according to the power-order.
It then chooses $(d,b)$ since $c$ dominates $d$ in the tree with both $(c,b)$ and $(d,b)$ added. 
}
\label{fig:example}
\end{figure}

We end this section by using the relation between parking function and orientations
to prove a property about the poset formed by parking functions. Recall given two parking functions
$f$ and $g$, we say that $f \prec g$ iff $f(v) \le g(v)$ for all vertices $v\in V$. Also recall 
the parking function $f \wedge g$ defined as $(f\wedge g)(v) := \min(f(v),g(v))$. Finally, recall ${\rm dom}(f) := \{g \mbox{ a parking function }: g\prec f\}$.

\begin{lem}
Given a non-maximum parking function $g$, let $F_g$ be the set of maximum parking functions that dominate $g$. Then,
$g = \bigwedge_{f\in F_g} f$.
\end{lem}
\begin{proof}
It is clear that any non-maximum parking function $g$ is dominated by the meet of all maximum parking functions which dominate it. 
The lemma claims that it is in fact exactly  equal to the meet. To show this, it suffices to show that for any $v$, there exists a parking function $f$ 
such that $f(v) = g(v)$, $f(u) \ge g(u)$ for all $u\in V$ and $f(w) > g(w)$ for exactly one vertex $w\neq v$. Continuing this process, we
get a maximum parking function which agrees with $g$ on $v$. Since this is true for all $v$, we are done.

To prove the above, note that from the bijection between parking functions and spanning trees, we see that any parking function $g$
uniquely corresponds to a {\em mixed} orientation of the edges of $G$, where a {\em mixed} orientation has some edges directed 
and others undirected. Moreover, the directed edges induce a DAG. 

Consider the mixed orientation with respect to $g$. Let $D$ be the DAG induced by the set of directed edges.
Suppose there is an edge $(u,w)$ which is not oriented where $u$ and $w$ are distinct from $v$; then in $D$ if there is a directed
path from $u$ to $w$, orient the edge from $u$ to $w$, else orient the edge from $w$ to $u$. if there is no path either way, orient in any of the two ways. We still have a mixed orientation where the directed edges are acyclic. Moreover, 
this increases the value of parking function on exactly one of $u$ or $w$, and keeps the value the same everywhere else.

Thus, the non-trivial case is when the only undirected edges are incident on $v$. Let $(u,v)$ be such an edge. Note that we do not
want to increase $g(v)$, that is, we want to orient $(u,v)$ without increasing $v$'s out-degree. In other words, we want to orient
it from $u$ to $v$. This is a problem if and only if there is a path from $v$ to $u$ in $D$. If so, consider  the longest such a path  $P = (v,w, \ldots ,u)$, where  $w$ is the neighbor of $v$ on $P$.

Note that apart from the edge $(v,w)$, there is no path from $v$ to $w$. If there were such a path, it can not use vertices from $P$, for $w$ ``dominates" (in the DAG) all vertices in $P$ other than $v$. Thus, the path from $v$ to $w$ must use  ``new" vertices making $P$ longer, and thus providing a contradiction.

Therefore, we can flip $(v,w)$ to $(w,v)$ and not create any cycles. That is,  the operation to get $f$ is to  flip $(v,w)$ to $(w,v)$ and orient $(v,u)$ from $v$ to $u$.
This only increases $g(w)$ but keeps everything else the same.
\end{proof}

\section{Product Graphs and $Q_n$-Parking Functions}

\setcounter{equation}{0}

Given two graphs $G_1$ and $G_2$, there is a standard notion of the Cartesian
product $G_1 \Box G_2$ of the two graphs.  Given a $G_1$-parking function and
a $G_2$-parking function, we define below a $G_1\Box G_2$-parking function
in a natural way that is symmetric in $G_1$ and $G_2$. 

\begin{defi}
Given $G_1 =(V_1, E_1)$ and $G_2=(V_2, E_2)$, the Cartesian product
graph $G_1\Box G_2 = (V,E)$
is  defined on the vertex set $V=V_1\times V_2$, using the edge set $E = E' \cup E''$, where
$E' = \bigl\{ \{(u_1,v), (u_2, v)\} : \{u_1,u_2\} \in E_1 \mbox{ and } v \in V_2 \bigr\}$\,, and
$E'' = \bigl\{ \{(u,v_1), (u, v_2)\} : \{v_1,v_2\} \in E_2 \mbox{ and } u \in V_1\bigr\}\,$.
\end{defi}
It is easy to see from the definition  that  the number of vertices in $G_1\Box G_2$ is  $|V_1|\,|V_2|$, and that  the number  of edges is 
 $|V(G_2)| \,|E(G_1)| + |V(G_1)| \,|E(G_2)|$.

\subsection{Parking functions on Product graphs}
\label{sec:PROD}
While it seems hard to characterize parking functions on $G_1 \Box G_2$, by simply
knowing those on $G_1$ and on $G_2$, the following result gives an explicit way to
construct a particular family of maximum parking functions on product graphs.

\begin{prop} 
\label{prop:PROD}
 Let $f_1 \in {\mathcal P}(G_1;q_1)$ and 
$f_2 \in {\mathcal P}(G_2;q_2)$. Then $f_1 \Box f_2 \in {\mathcal P}\bigl(G_1 \Box G_2; (q_1,q_2)\bigr)$, where $f_1\Box f_2(u,v) = f_1(u) + f_2(v) +1$, whenever $u\in G_1$
and $v\in G_2$.
 
 Further, if $f_1$ and $f_2$ are maximum parking functions, then $f_1 \Box f_2$ is a maximum parking function.\end{prop}

\Proof  Let $f=f_1\Box f_2$,  $G=G_1\Box G_2$, and $q=(q_1, q_2)$.
 To show that $f \in \PG$, once again we make crucial use of Dhar's marking algorithm; in particular, we will make use of 
the order in which the vertices of $G_1$ (and $G_2$) are marked in verifying
that   $f_1 \in \PGone$ (and $f_2\in \PGtwo$, respectively).
Using these in turn, we order the vertices in $G_1 \Box G_2$:
let $(u',v') <_\Box (u,v)$, if $u'$ is marked before $u$ in $G_1$, or if $u=u'$ and 
$v'$ is marked before $v$ in $G_2$.  We now prove that $f$ is a valid parking function, by showing that the vertices of $G_1\Box G_2$ can be marked, using Dhar, in precisely
the order given by $<_\Box$. 
We begin by noting that $f(q)= f(q_1, q_2) = f_1(q_1) +f_2(q_2) + 1 = -1$. Now consider
the vertices of $G$ inductively, using the order given by $<_\Box$. By the time the vertex $(u,v)$ is considered, observe that there are at least $f(u)+f(v)+2$ neighbors of $(u,v)$
that have already been marked, since they precede $(u,v)$ in $<_\Box$: indeed, at least $f(u)+1$ neighbors of the form
$(u',v)$ have been marked,  (since the graph induced by fixing the second coordinate $v$ is simply $G_1$), and similarly at least $f(v)+1$ neighbors of the form $(u,v')$ have also been marked.
Thus $(u,v)$ can be assigned the value $f(u,v) = f(u) + f(v) + 1$, and hence $f$ is a $G$-parking function.

Now, suppose that $f_1$ and $f_2$ are {\em maximum} $G_1$- and $G_2$-parking functions respectively.  Then, we must show that $f$ achieves the genus of ${\bf g}(G)$ minus one. This is easy to verify:
\begin{eqnarray*}
\|f_1 \Box f_2\| & = &  |V(G_2)|\, \|f_1\| +  |V(G_1)|\, \|f_2\| + |V(G_1)|\, |V(G_2)| \\
& = & |V(G_2)| \,|E(G_1)| + |V(G_1)| \,|E(G_2)| -  |V(G_1)|\, |V(G_2)| \\
& = & |E(G_1\Box G_2)| - |V(G_1\Box G_2)| = {\bf g}(G_1 \Box G_2) - 1. \vspace{-20mm}
\end{eqnarray*}

\subsection{$Q_n$-Parking Functions}
\label{sec:CUBE}

A quintessential product graph is the $n$-dimensional cube $Q_n$, obtained
by taking the product of an edge $Q_1$ with itself $n$ times.
For the purpose of this section, for integer $n\ge 1$,  we also view the $n$-cube $Q_n$ as the graph on $2^n$ vertices, which may conveniently be labeled by the $2^n$
binary vectors of length $n$, and with edges between vertices whose Hamming distance
is one. We are interested in understanding the parking functions on $Q_n$ with
respect to the vertex $q=(0,0, \ldots, 0)$.

\begin{defi}
\label{defn:can} For $n=1$, the unique parking function $f=f^1$ on $Q_1$ is canonical. 
For $n>1$, the parking function  
 $f^n = \underbrace{f \Box f \Box \cdots \Box f}_{n \text{ times}}$, obtained using the product graph construction, is defined as the {\em canonical} $Q_n$-parking function.
Further, if a $Q_n$-parking function $g$ is such that there exists a bijection $\phi:V(Q_n) \ra V(Q_n)$ such that $f(v)=g(\phi(v))$ for all $v \in V(Q_n)$, then we say that $g$ is {\em semi-canonical}.  
\end{defi}

  


Note that by Proposition~\ref{prop:PROD}, the canonical $Q_n$-parking function is a maximum parking function, and hence the semi-canonical one is also a maximum. 
   
\begin{exam}
{\em Not} all maximal $Q_n$-parking functions are semi-canonical.  For instance, consider the $Q_3$-parking function $f$ such that $f(000)=-1, f(001)=1, f(010)=0, f(100)=2, f(011)=0, f(101)=0, f(110)=0, f(111)=2$.  It is easy to verify that $f$ is a $Q_3$-parking function, but $f$ is not semi-canonical.  Since, $\|f\| = 4=|E|-|V|$,  $f$ is maximum (and thus maximal).
\end{exam}

Recall the partial order on parking functions, and the notion
of $\dom (f) = \{ g \in \PG : g \prec f\}$, for a parking function $f\in \PG$.

\begin{theo} 
\label{thm:domnum}
Let $f^n$ denote the canonical $Q_n$-parking function. Then $f^n(v)=\wgt(v)-1$,
where $\wgt(v)$ is the Hamming weight (the number of 1's in the 
binary representation) of the vertex  $v \in V(Q_n)$.
Consequently, if $f$ is semi-canonical, then $$ \left | \dom (f) \right | = \prod_{k=2}^n k^{n \choose k}.$$\end{theo}

\Proof The proof is by an easy induction on $n \in \Z^+$.
The base case is clear, since $f^1$ takes the values $-1, 0$.
For the induction step, for $n\ge1$,  write $f^{n+1} = f^n \Box f^1$, which 
by Proposition~\ref{prop:PROD} implies that, for $u\in V(Q_{n})$ and $v\in V(Q_1)=\{0,1\}$,
\[f^{n+1}(u,v)= f^n(u) + f^1(v) + 1\,.\]
Thus for $0\le k \le n$, the value $k$ can be obtained in $f^{n+1}$ either by taking a $k$ from $f^n$
and adding a zero to it (through, $f^1(0)+1 = -1+1$), or by taking a $k-1$ from
$f^n$ and adding a one to it (through, $f^1(1)=0+1$). By the induction hypothesis,
the number of $k$'s is
$${n \choose k} + {n \choose k-1}= {n+1 \choose k}\,,$$
completing the induction step. Also note that, by  definition, a semi-canonical
parking function also has the same distribution of integers.

To complete the proof of the theorem, recall  by Proposition~\ref{prop:dom_park}, 
that for any $v \in V(Q_n) \backslash \{q\}$, if $f(v)= k$ for $k \in \Z_{\geq 1}$, then for  $g\prec f$,  we may have $g(v)=0, 1, \ldots k$; this implies that there are $k+1$ possible values for such a $g$ with respect to $v$. 
Since $0 \leq k \leq n$ and each value in the range of $f$ is independent with respect to obtaining a dominated function $g$, we have that $$\left | \dom (f) \right | = \prod_{k=2}^n k^{n \choose k}.$$\eproof

\begin{rmk}
\label{can_semi} Note that it can directly be seen that $f$ is a $Q_n$-parking function, when $f$ is defined as $f(v)=\wgt(v)-1$.  Using Dhar's algorithm, starting
with $q =(0,0, \ldots , 0)$, we may proceed by marking vertices level by level (in the Boolean lattice
): Every vertex $v$ in level $k$ has precisely $k$ down-neighbors in level $k-1$, and they have all been marked, so $v$ can be marked and given value $k-1$.
Thus there are precisely ${n \choose k}$ vertices which obtain the value $k-1$, in such a canonical parking function on $Q_n$.
\end{rmk}
\begin{exam}
Note that {\em not} every semi-canonical parking function is canonical --  in satisfying the level by level property described in Remark~\ref{can_semi}. For instance, consider the $Q_3$-parking function $f$ such that $f(000)=-1, \  f(010)=f(100)=f(101)=0, \ f(001)=f(011)=f(111)=1$, and $f(110)=2$.  It is easy to check that $f$ is a $Q_3$-parking function, semi-canonical, but not canonical -- since, for example, the value 2 is adjacent to a 0.
\end{exam}

We hope the above remarks and examples indicate the difficulty in 
understanding the {\em maximum} parking functions on even a highly structured, 
symmetric graph such as the $n$-cube.  As far as we know, the number of maximum parking functions
of $Q_n$ is known only for $n\le 4$. For $n=2, 3$, and 4, this number is 3, 133, and 3040575, respectively.

\section{Diffuse states and acyclic orientations}
\label{sec:diffuse}
As mentioned in the introduction, in the context of chip-firing, the following
notion was introduced in \cite{JS95}. 

  \begin{defi}
 \label{defn:diffuse}
Given a connected graph $G$, a function $s:G\to  \Z^+$ is called a 
diffuse state if for every induced subgraph $G[A]=(A,E(A))
\subseteq G$, there exists some $u \in A$ such that $\deg(u)|_A \leq s(u)$.
Further, let $\|s\| := \sum_v s(v)$.
 \end{defi}

Note the (complementary) similarity with the definition of a parking function, 
by observing that $\deg(u)|_A \leq s(u)$ is equivalent to $\deg(u) - s(u) \leq \deg(u)|_{A^c}$. However, we have no special vertex such as $q$.
For chip-firing purposes, $s(v)$ may be thought of as the number of chips on $v$,
thus $\|s\|$ denotes the total number of chips in the graph.

First consider the following algorithm which constructs an acyclic orientation,
by using a given diffuse state $s$ with $|E(G)|$ chips. 

Step 1.  Since  $s$ is diffuse, we may find a vertex $v$ such
that $\deg(v)|_G = \deg(v) \leq s(v)$.

 Step 2.  Orient all of the edges incident to $v$ outward; delete $v$ and its incident edges.

 Step 3.  The resulting graph is diffuse since it is a subgraph of $G$; so we may
repeat Steps 1-2 until all edges of the graph are oriented.

Note that this process gives an acyclic orientation since we cannot orient
edges into a vertex which has out edges since this vertex has been deleted from
the graph. 

\begin{lem}
\label{lem:sinkset}
If $s$ is a diffuse state on graph $G$ with $\|s\|=|E(G)|$, then there is a vertex $v$ with  $s(v)=0$. For every diffuse state $s$,  the set $\{v: s(v)=0\}$ is an independent set in $G$.
\end{lem}

\proof
 The above algorithm which repeatedly removes vertices, removes at least
as many chips as the edges at each step. Before the last vertex, all edges
(hence all $|E(G)|$ chips) must have been removed, which means that the last vertex can has zero chips.
For the second part, if  $s(u)=s(v)=0$, and $u,v \in E(G)$ then the
set  $A=\{u,v\}$ violates the diffuse property.
\eproof


With a similar proof, it can also be shown that there exists a vertex $w$ 
such that  $s(w)=\deg(w)$, under the hypothesis of the above lemma.

\begin{theo}
\label{thm:diffuse}
There is a bijection between the set ${\rm  D}(G)$ of diffuse states with $|E(G)|$ chips and the set ${\rm AO}(G)$ of acyclic orientations of a connected graph $G$.
\end{theo}

\proof
Let $E(G) = m$. 
The proof is based on two injections between the sets, going in each direction.
First, given an acyclic orientation ${\mathcal O} \in {\rm AO}(G)$, define the nonnegative function $s=s_{\mathcal O}$ by letting $s(v)$ be the out-degree of $v$ in the orientation ${\mathcal O}$. Clearly, the mapping is
one-one, since the out-degree sequence uniquely determines an acyclic orientation; also  $\|s\|=m$, since the sum of out-degrees equals the number of edges.  
To see
that $s$ is diffuse, simply observe that,  each induced subgraph
$G[A] \subseteq G$ has a (local) source $y \in A$ when  restricted to the acyclic orientation induced on $A$; such a source  $y$ satisfies
 $s(y) \geq \deg(y)|_A$, since the out-degree of  $y$ is at least the degree $\deg(y)|_A$.

For an injection in the other direction, we make use Lemma~\ref{lem:sinkset}.
Given a diffuse state $s$ with $m$ chips,
we construct an acyclic orientation, by constructing a DAG:
 Lemma~\ref{lem:sinkset} guarantees the existence of sink(s); so we construct the orientation, by (i) repeatedly removing the current set
of sinks, and (ii) subtracting a chip from each in-neighbor of a removed sink.
It is easy to see that the updated function $s$ at each step is still a diffuse
state on the remaining graph.  Note that this construction is one-one:  for $s, s'$ different diffuse states, 
simply consider the first time the current sets of sinks (in the DAGs) differ, when we start with $s$ versus $s'$; since the underlying graph is the same, there must be such a time whenever $s\neq s'$.
\eproof

\begin{rmk}
Thanks again to an anonymous referee of \cite{BT08}, we learnt the following: that diffuse states are also  in 1-1 correspondence with the critical configurations of the chip-firing game of Bj\"orner-Lov\'asz-Shor \cite{BLS91}, and that the first part of the above proof appears as Theorem~3.3, part (b) in \cite{BLS91}. In addition, it is easy  to go between the above theorem and Theorem~\ref{thm:biject}: simply add a new vertex $q$ adjacent to every vertex in $G$ and obtain a new graph 
$\hat{G}$; then $\hat{G}$-parking functions (with respect to $q$, say) correspond to diffuse states in $G$, and acyclic orientations with a unique source at $q$ in $\hat{G}$ correspond to acyclic orientations of $G$.
\end{rmk}

\section{Concluding Remarks}
\label{sec:CONCL}
Soon after the completion of this work, Igor Pak kindly pointed us to the work of Olivier Bernardi \cite{B08},
where bijective proofs are derived for interpretations of each of the evaluations of the Tutte polynomial ${\mathcal T}_G(x,y)$, for $0\le x, y \le 2$, in terms of orientations. A key to this seems to be a nice combinatorial embedding of (the edges of) the  graph and a rewriting of the Tutte polynomial using notions of internal and external {\em embedding} activity. 

In addition to the questions mentioned in the previous sections, several challenging
problems remain open.  Given an arbitrary graph $G$, it is a classical open problem in the topic of Markov chain Monte Carlo (MCMC) algorithms \cite{JS96}, to efficiently  generate
an acyclic orientation {\em uniformly at random} from the set of all such orientations.
Due to the observations above, a closely related problem would be
to generate at random an acyclic orientation with a uniquely identified  sink (or source), or equivalently,  to sample uniformly from the set of {\em safe} spanning trees of $G$. 
The MCMC technique suggests the following natural approach to this problem:
it is well known (see for example, \cite{FM92,FR05}) that the so-called bases exchange walk provides an efficient way to sample uniformly from the set of {\em all} spanning trees of a given graph $G$.  However it remains to be seen whether (and how)  restricting such a random walk to the set of safe trees affects the {\em mixing time} of the walk --
 the time by which the walk converges to its steady state distribution,  uniform 
on the set of safe trees. 
Given that the exact enumeration of the number of safe trees of $Q_n$ is also open,
an interesting first step might be to analyze such a walk on the trees of $Q_n$.

Independent of the above approach, other ways of providing asymptotically accurate estimates, for large $n$, of  the number (or even the logarithm of the number) of maximum parking functions on $Q_n$\,,  remains  interesting and presumably a challenging exercise.

\noindent {\bf Acknowledgements.}  The authors thank Matt Baker for introducing them to  $G$-parking functions and for helpful discussions.  The authors also  thank Adam Marcus for helpful discussions, and Herb Wilf for providing the chromatic polynomial of $Q_4$, which identified the number of maximum parking functions of $Q_4$.  The authors are grateful to the anonymous referees for making several critical remarks which helped improve the presentation here. The present work originated from a summer 2007 REU, supported by an REU supplement to the last author's NSF grant DMS-0701043.

\end{document}